\documentclass[12pt]{amsart} %amsfonts,
\usepackage{amsmath,amssymb,amscd,amsthm}
\usepackage{latexsym}
\usepackage{graphicx}
\usepackage[english]{babel}
\usepackage[latin1]{inputenc}       % allow Latin1 characters

\usepackage{textcomp}

\usepackage{pgf,pgfarrows,pgfnodes,pgfautomata,pgfheaps}
\usepackage{colortbl}

\newtheorem{theorem}{Theorem}[section]

\newtheorem{lemma}[theorem]{Lemma}

\newtheorem{defn}[theorem]{Definition}

\newtheorem{remark}[theorem]{Remark}

\def\irr#1{{\rm Irr}(#1)}
\def\irrr#1#2 {\irr {#1 \mid #2}}
\newcommand{\R}{\mathbb R}

\newcommand{\F}{e^{-\varphi(|y|)}}
\newcommand{\E}{\textbf{E}}
\newcommand{\V}{\textbf{S}}

\begin{document}

\title[Surface area of polytopes]{Maximal surface area of polytopes with respect to log-concave rotation invariant measures.}
  \keywords{ convex bodies, convex polytopes, Surface area, Gaussian measures\\
   Supported in part by U.S. National Science Foundation Grant DMS-1101636.}
\author{Galyna Livshyts}
%\address{Department of Mathematics, Kent State University,
%Kent, OH 44242, USA} \email{glivshyt@kent.edu}
\date{}

\begin{abstract}
It was shown in \cite{GL} that the maximal surface area of a convex set in $\R^n$ with respect to a rotation invariant log-concave probability measure $\gamma$ is of order $\frac{\sqrt{n}}{\sqrt[4]{Var|X|} \sqrt{\mathbb{E}|X|}}$, where $X$ is a random vector in $\R^n$ distributed with respect to $\gamma$. In the present paper we discuss surface area of convex polytopes $P_K$ with $K$ facets. We find tight bounds on the maximal surface area of $P_K$ in terms of $K$. We show that $\gamma(\partial P_K)\lesssim \frac{\sqrt{n}}{\mathbb{E}|X|}\cdot\sqrt{\log K}\cdot\log n$ for all $K$. This bound is better then the general bound for all $K\in [2,e^{\frac{c}{{\sqrt{Var|X|}}}}]$. Moreover, for all $K$ in that range the bound is exact up to a factor of $\log n$: for each $K\in [2,e^{\frac{c}{{\sqrt{Var|X|}}}}]$ there exists a polytope $P_K$ with at most $K$ facets such that $\gamma(\partial P_K)\gtrsim \frac{\sqrt{n}}{\mathbb{E}|X|}\sqrt{\log K}.$ %For the measures $\gamma_p$ with densities $C_{n,p} e^{-\frac{|y|^p}{p}}$ (where $p>0$) we obtain: $\gamma_p(\partial P_K)\lesssim \frac{\sqrt{n}}{\mathbb{E}X}\sqrt{\log K},$ which was obtained for the standard Gaussian measure $\gamma_2$ by F. Nazarov.
\end{abstract}
\maketitle

\section{Introduction}

In this paper we study properties of the surface area of convex polytopes with respect to log-concave rotation invariant probability measures. For sets $A, B\subset\R^n$ the Minkowski sum is defined as
$$A+B=\{a+b\,|\,a\in A, b\in B\}.$$

For a scalar $\lambda$ the dilated set is
$$\lambda A:=\{\lambda a\,|\, a\in A\}.$$

A measure $\gamma$ on $\R^n$ is called log-concave if for any measurable sets $A, B\subset\R^n$ and for any $\lambda\in [0,1]$,
$$\gamma(\lambda A+(1-\lambda)B)\geq \gamma(A)^{\lambda}\gamma(B)^{1-\lambda}.$$
It was shown by Borrell \cite{bor}, that a measure is log-concave if and only if it has a density with respect to the Lebesgue measure on some affine hyperplane, and this density is a log-concave function.  Log-concave measures have been studied intensively in the recent years. For the background and numerous interesting properties, see for example \cite{Kl}, \cite{Kl2}, \cite{KM} and \cite{MP}.

A measure $\gamma$ is called rotation invariant if, for every rotation $T$ and for every measurable set $A,$
$$\gamma(TA)=\gamma(A).$$
Log-concave rotation invariant measures appear for example in \cite{KM}, \cite{bob1}, \cite{bob2}, \cite{bob3} and \cite{GL}.

In the present paper we restrict our attention to probability measures (that means that the measure of the whole space is equal to $1$). Examples of log-concave rotation invariant probability measures are the Standard Gaussian Measure $\gamma_2$ and the Lebesgue measure restricted on a ball.

Let $X$ be a random vector in $\R^n$ distributed with respect to a measure $\gamma$. We introduce
\begin{equation}\label{E}
\E:=\mathbb{E}|X|
\end{equation}
and
\begin{equation}\label{V}
\V:=\frac{\sqrt{\mathbb{E}\left(|X|-\mathbb{E}|X|\right)^2}}{\mathbb{E}|X|},
\end{equation}
the expectation and the normalized standard deviation of the absolute value of $X$. $\E$ and $\V$ are natural parameters of the measure $\gamma$. For rotation invariant measures $\V\in [\frac{c_1}{n},\frac{c_2}{\sqrt{n}}]$, where $c$ and $c'$ are absolute constants (see \cite{Kl} or \cite{GL}, Remark 2.9). The parameter $\V$ is closely related to $\sigma=\sqrt{\mathbb{E}\left(|X|-\mathbb{E}|X|\right)^2}$. The famous Thin Shell Conjecture suggests that $\sigma$ is bounded from above by an absolute constant for all isotropic (see, for example, \cite{MP}, \cite{Kl} for definitions and properties) log-concave measures (see \cite{Kl}, \cite{Kl2}, \cite{Kl3}, \cite{Fl1}, \cite{Fl2}, \cite{Ron}, \cite{GUMI}). Currently, the best bound is $Cn^{\frac{1}{3}}$ and is due to Gudeon and E. Milman \cite{GUMI}.

The Minkowski surface area of a convex set $Q$ with respect to the measure $\gamma$ is defined to be
\begin{equation}\label{outer}
\gamma(\partial Q)=\liminf_{\epsilon\rightarrow +0}\frac{\gamma((Q+\epsilon B_2^n)\backslash Q)}{\epsilon},
\end{equation}
where $B_2^n$ denotes Euclidian ball in $\R^n.$ In many cases the Minkowski surface area has an integral representation:
\begin{equation}\label{mink_sa_integral}
\gamma(\partial Q)=\int_{\partial Q} f(y) d\sigma(y),
\end{equation}
where $f(y)$ is the density $\gamma$ and $d\sigma(y)$ stands for the Lebesgue surface measure on $\partial Q$ (see for example Appendix of \cite{kane} and Appendix of \cite{GL}).

The questions of estimating the surface area of $n-$dimensional convex sets with respect to the Standard Gaussian Measure have been actively studied. Sudakov, Tsirelson \cite{ST} and Borell \cite{B} proved, that among all convex sets of a fixed Gaussian volume, half spaces have the smallest Gaussian surface area. Mushtari and Kwapien asked the reverse version of the isoperimetric inequality, i.e. how large the Gaussian surface area of a convex set $Q\subset \R^n$ can be. It was shown by Ball \cite{ball}, that Gaussian surface area of a convex set in $\R^n$ is asymptotically bounded by $C n^{\frac{1}{4}}$, where $C$ is an absolute constant. Nazarov \cite{fedia} proved the sharpness of Ball's result and gave the complete solution to this asymptotic problem:
\begin{equation}\label{fedia}
0.28 n^{\frac{1}{4}}\leq \max_{Q\in \mathcal{K}_n} \gamma_2(\partial Q)\leq 0.64 n^{\frac{1}{4}},
\end{equation}
where by $\mathcal{K}_n$ we denote the set of all convex sets in $\R^n$.

Further estimates for $\gamma_2(\partial Q)$ for the special case of polynomial level set surfaces were provided by Kane \cite{kane}. He showed that for any polynomial $P(y)$ of degree $d$, $\gamma_2(P(y)=0)\leq \frac{d}{\sqrt{2}}$.

For the case of all rotation invariant log-concave measures it was shown in \cite{GL}, that
$$\max_{Q\in \mathcal{K}_n}\gamma(\partial Q)\approx \frac{\sqrt{n}}{\E \cdot\sqrt{\V}}.$$

Let $K$ be a given positive integer. In the present paper we consider the family of $n-$dimensional convex polytopes with $K$ facets, where by a ``polytope'' we mean the intersection of $K$ half-spaces (we do not assume compactness as it is irrelevant for the type of questions we consider). We obtain the bounds on the surface area of the polytope with $K$ facets with respect to rotation invariant log-concave measure $\gamma$ in terms of $K$ and the natural parameters of $\gamma$.

It is not hard to show that the $\gamma-$surface area of any half-space does not exceed $C\frac{\sqrt{n}}{\E}$, for some absolute constant $C$ (see (\ref{hyperplane}) below). Thus the immediate bound for the surface area of a polytope with $K$ facets is $C\frac{K\sqrt{n}}{\E}$. In the present paper we show a sharper estimate from above. We also show an estimate from below on the maximal surface area of a convex polytope with $K$ facets. Both of the estimates match up to a $\log n$ factor.

The estimate from above is the content of the following Theorem:

\begin{theorem}\label{1}
Let $n\geq 2.$ Fix positive integer $K\in[2,e^{\frac{c}{\V}}]$. Let $P$ be a convex polytope in $\R^n$ with at most $K$ facets. Let $\gamma$ be a rotation invariant log-concave measure with $\E$ and $\V$ defined by (\ref{E}) and (\ref{V}). Then
$$\gamma(\partial P)\leq C\frac{\sqrt{n}}{\E}\cdot\sqrt{\log K}\cdot\log\frac{1}{\V \log K},$$
where $C$ and $c$ stand for absolute constants.
\end{theorem}

\begin{remark}
We note that for $K\geq e^{\frac{c}{\V}}$ the bound from Theorem \ref{1} becomes worse then the general bound $\frac{\sqrt{n}}{\sqrt[4]{Var|X|} \sqrt{\mathbb{E}|X|}}$. The latter bound is also optimal for all $K\geq e^{\frac{c}{\V}}$.
\end{remark}

\begin{remark}
Since $\V\in [\frac{c}{n},\frac{C}{\sqrt{n}}]$ for all log-concave rotation invariant measures, Theorem \ref{1} reads in fact, that for all $K\in[1,e^{\frac{c}{\V}}]$,
$$\gamma(\partial P)\leq C\frac{\sqrt{n}}{\E}\cdot\sqrt{\log K}\cdot\log\frac{n}{\log K}\leq C\frac{\sqrt{n}}{\E}\cdot\sqrt{\log K}\cdot\log n.$$
\end{remark}

Theorem \ref{1}, up to a log factor, is a generalization of the following Theorem of Nazarov \cite{PK}:

\begin{theorem}[\textbf{F. Nazarov}]\label{th_naz}
Let $n\geq 2$ and $K\geq 2$ be integers. Let $P$ be a convex polytope in $\R^n$ with at most $K$ facets. Let $\gamma_2$ be the Standard Gaussian Measure. Then there exist a positive constant $C$ such that
$$\gamma_2(\partial P)\leq C\sqrt{\log K}.$$
\end{theorem}
For a generalization of the Theorem of Nazarov in an entirely different set up see \cite{kane2}. See also Section 5 of the present paper for the proof of the analogous result for measures with densities $C(n,p)e^{-\frac{|y|^p}{p}}$. The case $p=2$ corresponds to the Gaussian measure. Theorem \ref{final} from Section 5 is a generalization of the Theorem of Nazarov.

We also obtain a lower bound for the maximal surface area of a convex polytope with $K$ facets. It proves sharpness of Theorem \ref{th_naz} of Nazarov. It also shows sharpness of Theorem \ref{1} up to a $\log n$ factor:

\begin{theorem}\label{2}
Let $n\geq 2.$ Let $\gamma$ be a rotation invariant log-concave measure with $\E$ and $\V$ defined by (\ref{E}) and (\ref{V}). Fix positive integer $K\in[2,e^{\frac{c}{\V}}]$. Then there exists a convex polytope $P$ in $\R^n$ with at most $K$ facets such that
$$\gamma(\partial P)\geq C'\frac{\sqrt{n}}{\E}\sqrt{\log K},$$
where $c$ and $C'$ stand for absolute constants.
\end{theorem}

The next section is dedicated to some technical preliminaries. In Section 3 we give the proof of Theorem \ref{1}. In Section 4 we prove Theorem \ref{2}. Finally, in Section 5 we show that Theorem \ref{1} can be refined in some partial cases of measures, which include the Standard Gaussian measure.

\section{Preliminaries and definitions}

This section is dedicated to some general properties of rotation invariant log-concave measures.
We outline some elementary facts which are needed for the proof. Some of them have appeared in literature.
See \cite{Kl} for an excellent overview of the properties of log-concave measures; see also \cite{GL} for more details and the proofs of the facts listed in the present section.

We use notation $|\cdot|$ for the norm in Euclidean space $\R^{n}$; $|A|$ stands for the Lebesgue measure  of a measurable set $A \subset \R^{n}$.
We write     $B_2^{n}=\{x\in \R^{n}: |x| \le 1\}$ for the unit ball in $\R^{n}$ and $\mathbb{S}^{n-1}=\{x\in \R^{n}: |x|=1\}$ for the unit sphere. We denote $\nu_{n}=|B_2^{n}|=\frac{\pi^{\frac{n}{2}}}{\Gamma(n/2+1)}$.

We shall use notation $\precsim$ for an asymptotic inequality: we say that $A(n)\precsim B(n)$ if there exists an absolute positive constant $C$ (independent of $n$), such that $A(n)\leq C\cdot B(n)$. Correspondingly, $A(n)\approx B(n)$ means that $B(n)\precsim A(n)\precsim B(n)$. Also in the present paper $C$, $c,$ $c_1$ etc denote absolute constants which may change from line to line.

We fix a convex nondecreasing function $\varphi(t):[0,\infty)\rightarrow [0,\infty]$. Let $\gamma$ be a probability measure on $\R^{n}$ with density $C_{n}\F$. The normalizing constant $C_{n}$ equals to $[n\nu_{n} J_{n-1}]^{-1}$, where
\begin{equation}\label{maint}
J_{n-1}=\int_0^{\infty} t^{n-1} \F dt.
\end{equation}
The measure $\gamma$ is rotation invariant and log-concave; conversely, every rotation invariant log-concave measure is representable this way in terms of some convex function $\varphi$. Since we normalize the measure anyway, we may assume that $\varphi(0)=0$. We will be aiming for the estimates for
$$\gamma(\partial P)=[n\nu_{n} J_{n-1}]^{-1}\int_{\partial P} e^{-\varphi(|y|)} d\sigma(y),$$
where $P$ is a convex polytope with $K$ facets.

Without loss of generality we assume that $\varphi\in C^2[0,\infty)$. This can be shown by the standard smoothing argument (see, for example, \cite{cord}).

We introduce the notation
\begin{equation}\label{g}
g_{n-1}(t)=t^{n-1} e^{-\varphi(t)}.
\end{equation}
\begin{defn}\label{tzerodef}
We define $t_0$ to be the point of maxima of the function $g_{n-1}(t)$, i.e., $t_0$ is the solution of the equation
\begin{equation}\label{to}
\varphi'(t)t=n-1.
\end{equation}
\end{defn}
The equation (\ref{to}) has a solution, since $t\varphi'(t)$ is non-decreasing, continuous and $\lim_{t\rightarrow +\infty} t\varphi'(t)=+\infty$.
This solution is unique, since $t\varphi'(t)$ strictly increases on its support. This definition appears in most of the literature dedicated to
spherically symmetric log-concave measures: see, for example, \cite{KM} or \cite{Kl} (Lemma 4.3), as well as \cite{GL} (Definition 2).

The following Lemma is proved in \cite{KM}. It provides asymptotic bounds for $J_{n-1}$.
\begin{lemma}\label{mainintegral}
$$\frac{g_{n-1}(t_0)t_0}{n} \leq J_{n-1}\leq \sqrt{2\pi}(1+o(1))\frac{g_{n-1}(t_0)t_0}{\sqrt{n-1}}.$$
\end{lemma}

%The function $g_n(t)=t^n e^{-\varphi(t)}$ is log-concave on $[0,\infty)$, and we shall apply Lemma \ref{log_concave} with $g(t)=g_n(t)$ and $\psi=1$.
The following definitions appear in \cite{GL} (Definition 3).
\begin{defn}\label{lambdas}
Define the ''outer'' $\lambda_o$ to be a positive number satisfying:
\begin{equation}\label{lambda_def_outer}
\varphi(t_0(1+\lambda_o))-\varphi(t_0)-(n-1)\log(1+\lambda_o)=1.
\end{equation}
Similarly, define the ''inner'' $\lambda_i$ as follows:
\begin{equation}\label{lambda_def_inner}
\varphi(t_0(1-\lambda_i))-\varphi(t_0)-(n-1)\log(1-\lambda_i)=1.
\end{equation}
We put
\begin{equation}\label{lambda}
\lambda:=\lambda_i+\lambda_o.
\end{equation}
\end{defn}

We note that (\ref{lambda_def_outer}) is equivalent to
\begin{equation}\label{lambda_def_outer1}
g_{n-1}(t_0)=e\cdot g_{n-1}(t_0(1+\lambda_o)),
\end{equation}
and (\ref{lambda_def_inner}) is equivalent to
\begin{equation}\label{lambda_def_inner1}
g_{n-1}(t_0)=e\cdot g_{n-1}(t_0(1-\lambda_i)).
\end{equation}

Parameter $\lambda$ from (\ref{lambda}) has a nice property.
\begin{lemma}\label{int_lam}
$$J_{n-1} \approx\lambda t_0 g_{n-1}(t_0).$$
\end{lemma}
See \cite{GL} (Lemma 4) for the details and the proof. The following fact is also presented in \cite{GL} (Lemma 5).
\begin{lemma}\label{moments}
For all $n\geq 2$,
$$\frac{J_{n}}{J_{n-1}}\approx t_0.$$
\end{lemma}

The above implies, that $t_0\approx \E=\mathbb{E}|X|$, where $X$ is a random vector in $\R^{n}$ distributed with respect to $\gamma.$ Also, $\lambda\approx \V,$ where $\V$ is defined by (\ref{V}) (see \cite{GL} (Lemmas 9 and 10) for the details).

\begin{remark}\label{remarklambda}
We note that Lemma \ref{int_lam} together with Lemma \ref{mainintegral} imply that $\lambda\in [\frac{c_1'}{n},\frac{c_2'}{\sqrt{n}}]$. Both of the estimates are exact: it is equal to $\frac{c}{n}$ for Lebesgue measure concentrated on a ball and to $\frac{C}{\sqrt{n}}$ for Standard Gaussian Measure.

Since $\lambda\approx\V,$ we claim also that $\V\in [\frac{c_1}{n},\frac{c_2}{\sqrt{n}}]$.
\end{remark}

We are now after the restated versions of the Theorems \ref{1} and \ref{2}:
\begin{theorem}\label{3}
Let $n\geq 2.$ Let $K\in [2, e^{\frac{c}{\lambda}}]$. Let $P$ be a convex polytope in $\R^n$ with at most $K$ facets. Let $\gamma$ be a rotation invariant log-concave measure. Then
$$\gamma(\partial P)\leq C\frac{\sqrt{n}}{t_0}\cdot\sqrt{\log K}\cdot\log\frac{1}{\lambda\log K}.$$
%where $C$ stands for an absolute constant.
\end{theorem}
and
\begin{theorem}\label{4}
Let $n\geq 2.$ Let $\gamma$ be a rotation invariant log-concave measure. Fix positive integer $K\in[2,e^{\frac{c}{\lambda}}]$. Then there exists a polytope $P$ in $\R^n$ with at most $K$ facets such that
$$\gamma(\partial P)\geq C'\frac{\sqrt{n}}{t_0}\sqrt{\log K}.$$
%where $C'$ stands for an absolute constant.
\end{theorem}

The following Lemma is an elementary fact about log-concave functions (for example, it appears in \cite{GL} as Lemma 3). %We will use it to estimate the ''tails'' of $J_n$.
\begin{lemma}\label{log_concave}
Let $g(t)=e^{f(t)}$ be a log-concave function on [a,b] (where both $a$ and $b$ may be infinite). We assume that $f\in C^2[a,b]$. Let $t_0$ be the point of maxima of $f(t)$. Assume that $t_0>0.$ Consider $x>0$ and $\psi>0$ such that
$$f(t_0)-f((1+x)t_0)\geq \psi.$$
Then,
$$\int_{(1+x)t_0}^{b} g(t)dt\leq \frac{x t_0 g(t_0)}{\psi e^{\psi}}.$$
Similarly, if $f(t_0)-f((1-x)t_0)\geq \psi,$
$$\int_a^{(1-x)t_0} g(t)dt\leq \frac{x t_0 g(t_0)}{\psi e^{\psi}}.$$
\end{lemma}

The next Lemma is similar to Lemma 12 from \cite{GL}.

\begin{lemma}\label{annulus}
Pick $\psi\in [1, c\log \frac{1}{\lambda}]$. Define $\mu$ to be smallest positive number such that
\begin{equation}\label{condition_lem_an}
\varphi\left(t_0(1+\mu)\right)-\varphi(t_0)-(n-1)\log(1+\mu)\geq\psi.
\end{equation}
%Assume that $\mu=o(n^{-\frac{1}{3}})$.
Define
$$A:=(1+\mu)t_0 B_2^{n}\setminus \frac{t_0}{2e}B_2^{n}.$$
%and denote the compliment to the set $A$ by $A^c$.
We claim, that such $\mu$ is well-defined and
$$\gamma(\partial Q\setminus A)\precsim \frac{\sqrt{n}}{t_0\lambda \sqrt{\psi} e^{\psi}}.$$
\end{lemma}
\noindent\textbf{Proof.} First, consider $M=Q\cap \frac{t_0}{2e}B_2^{n+1}$. Then,
$$\gamma(M)\leq \frac{1}{(n-1) \nu_{n} J_{n-1}} \int_{M} e^{-\varphi(|y|)} d\sigma(y)\leq \frac{|M|}{(n-1) \nu_{n} J_{n-1}}\leq$$
\begin{equation}\label{estimate1}
\frac{|\frac{t_0}{2e} \mathbb{S}^{n-1}|}{(n-1) \nu_{n} J_{n-1}}\approx \frac{t_0^{n-1}}{(2e)^{n-1} \lambda t_0 e^{-\varphi(t_0)}t_0^{n-1}}=\frac{1}{\lambda t_0}\cdot \frac{e^{\varphi(t_0)}}{(2e)^{n-1}},
\end{equation}
where the equivalency follows from Lemma \ref{int_lam} and (\ref{annulusrough}). By the Mean Value Theorem, $\varphi(t_0)\leq n$, so we estimate (\ref{estimate1})
from above by  $c\frac{2^{-n}}{\lambda t_0}$. In a view of Remark \ref{remarklambda}, the latter bound is much better then the one stated in the Lemma.

Next, let $N=\partial Q\setminus (1+\mu)t_0 B_2^{n}.$
%For $\mu=o(n^{-\frac{1}{3}})$,
In the current range of $\psi,$ we observe:
\begin{equation}\label{psimu}
\psi=\frac{(n-1)\mu^2}{2}+o(1).
\end{equation}

We obtain the following integral expression for $e^{-\varphi(|y|)}$ (inspired by \cite{ball}):
$$e^{-\varphi(|y|)}=\int_{|y|}^{\infty} \varphi'(t)e^{-\varphi(t)}dt=\int_0^{\infty} \varphi'(t)e^{-\varphi(t)}\chi_{[0,t]}(|y|)dt,$$
where $\chi_{[0,t]}$ stands for characteristic function of the interval $[0,t]$. In the current range of $|y|$,
$$e^{-\varphi(|y|)}=\int_{(1+\mu)t_0}^{\infty} \varphi'(t)e^{-\varphi(t)}\chi_{[0,t]}(|y|)dt.$$
Using the above, passing to the polar coordinates and integrating by parts, we get
$$\gamma(N)\leq \frac{1}{J_{n-1}}\int_{(1+\mu)t_0}^{\infty} t^{n-1} \varphi'(t)e^{-\varphi(t)}dt\approx$$
\begin{equation}\label{mu_comp}
\frac{g_{n-1}((1+\mu)t_0)+(n-1)\int_{(1+\mu)t_0}^{\infty} g_{n-2}(t)dt}{\lambda t_0 g_{n-1}(t_0)}.
\end{equation}
Lemma \ref{log_concave}, applied with $x=\mu$ and $\psi$, together with (\ref{psimu}) entails that (\ref{mu_comp}) is asymptotically less than
$$
\frac{e^{-\psi}}{\lambda t_0}+\frac{n\mu}{\lambda t_0\psi e^{\psi}}\leq \frac{1}{\lambda t_0} \cdot (1+\frac{\mu n}{\psi}) e^{-\psi}\precsim \frac{\sqrt{n}}{t_0\lambda \sqrt{\psi} e^{\psi}},
$$
which implies the estimate. $\square$

We use Lemma \ref{annulus} with $\mu\approx \frac{\sqrt{\log{\frac{1}{\lambda\sqrt{\log K}}}}}{\sqrt{n}}\leq 1$. We get, that for $$A:=\left(1+\frac{\sqrt{\log{\frac{1}{\lambda\sqrt{\log K}}}}}{\sqrt{n}}\right)t_0 B_2^{n}\setminus \frac{t_0}{2e}B_2^{n},$$
it holds that
$$\gamma(\partial Q\setminus A)\precsim \frac{\sqrt{n}}{t_0}\sqrt{\log K}.$$
%$$\gamma(\partial Q\setminus 5t_0 B_2^n)\lesssim \frac{e^{-n}}{t_0},$$
%and
%$$\gamma(\partial Q\cap \frac{t_0}{e^2} B_2^n)\lesssim \frac{2^{-n}}{t_0}.$$
Let $y\in\partial Q$. In an account of the above, we may assume that
\begin{equation}\label{annulusrough}
|y|\approx t_0
\end{equation}
throughout the proof.

We consider the hyperplane $H$ passing through the origin.
\begin{equation}\label{hyp}
\gamma(H)\approx\frac{1}{n\nu_{n}J_{n-1}}\int_{\R^{n-1}} e^{-\varphi(|y|)} d\sigma(y)=\frac{(n-1)\nu_{n-1}J_{n-2}}{n\nu_n J_{n-1}}.
\end{equation}

It is well known that
\begin{equation}\label{balls_ratio}
\frac{\nu_{n-1}}{\nu_{n}}\approx \sqrt{n}.
\end{equation}
%The simplest way to see it is to write a volume of $B_2^n$ by Fubbini's Theorem:
%$$\nu_{n+1}=|B_2^{n+1}|=\int_0^1 (1-t^2)^{\frac{n}{2}}\nu_{n}dt\approx\frac{\nu_n}{\sqrt{n}}.$$
Applying (\ref{balls_ratio}) together with Lemma \ref{moments} and (\ref{hyp}), we obtain that
\begin{equation}\label{hyperplane}
\gamma(H)\approx\frac{\sqrt{n}}{t_0}.
\end{equation}
Thus the trivial bound on the surface area of a polytope $P$ with $K$ facets in $\R^{n}$ is $\frac{\sqrt{n}}{t_0}K=\frac{\sqrt{n}}{\E}K$. We shall improve it.

\section{Proof of the upper bound part}

Let $Q$ be a convex set in $\R^n$. For $y\in\partial Q$ define
\begin{equation}\label{alpha}
\alpha(y):=\cos(y,n_y),
\end{equation}
where $n_y$ stands for the normal vector at $y.$ We also define
\begin{equation}\label{psi}
\psi(y):=\log\frac{g_{n-1}(t_0)}{g_{n-1}(|y|)}=\varphi(t_0)-\varphi(|y|)-(n-1)\log\frac{|y|}{t_0}.
\end{equation}

It was shown in \cite{GL} (Equation (46)) that

\begin{equation}\label{system1}
\gamma(\partial Q)\lesssim \max_{y\in \partial Q}\frac{1}{\lambda|y|\alpha(y)e^{\psi(y)}}.
\end{equation}

It was also shown in \cite{GL} (Equation (49) and Proposition 1) that

\begin{equation}\label{system2}
\gamma(\partial Q)\lesssim \max_{y\in \partial Q}\sqrt{n}\sqrt{\psi(y)}\frac{\alpha(y)\sqrt{n}\sqrt{\psi(y)}+1}{|y|}.
\end{equation}

Pick any $\theta\in \mathbb{S}^{n-1}$ and $\rho>0$. Let
$$H_{\rho}=\{x\in\R^n\,|\, \langle x, \theta\rangle=\rho\}$$
be a hyperplane at distance $\rho$ from the origin. We note that for $y\in H_\rho,$ $\alpha(y)|y|=\rho.$ So we introduce another function
\begin{equation}\label{r(y)}
r(y):=\frac{\sqrt{n}}{t_0}\alpha(y)|y|.
\end{equation}
For all $y\in H_{\rho}$ the function $r(y)=\frac{\sqrt{n}}{t_0}\rho$. Applying (\ref{annulusrough}), we rewrite (\ref{system1}) and (\ref{system2}) in terms of $r(y)$:

\begin{equation}\label{system1_r}
\gamma(\partial Q)\lesssim  \frac{\sqrt{n}}{t_0}\max_{y\in \partial Q}\frac{1}{\lambda r(y) e^{\psi(y)}},
\end{equation}

\begin{equation}\label{system2_r}
\gamma(\partial Q)\lesssim \frac{\sqrt{n}}{t_0}\max_{y\in \partial Q} \left(r(y)\psi(y)+\sqrt{\psi(y)}\right).
\end{equation}

We are going to estimate the measure of each facet using both (\ref{system1_r}) and (\ref{system2_r}), and ``the breaking point'' is going to depend on how far the facet is from the origin. So we minimize the expression
\begin{equation}\label{tominimize}
\frac{1}{\lambda r e^{\psi}}+r\psi+\sqrt{\psi}
\end{equation}
in $\psi$ in terms of $r.$ If $r\lesssim \frac{1}{\sqrt{\psi}}$, the minimum of (\ref{tominimize}) is equivalent to the minimum of $\frac{1}{\lambda r e^{\psi}}+\sqrt{\psi}$ which is achieved when $\psi\approx \log\frac{1}{\lambda r}$ and is approximately equal to $\sqrt{\log\frac{1}{\lambda r}}$. If $r\precsim \frac{1}{\sqrt{\psi}}$, the minimum of (\ref{tominimize}) is equivalent to the minimum of $\frac{1}{\lambda r e^{\psi}}+r\psi$ which is achieved when $\psi\approx \log\frac{1}{r^2\lambda}$ and is approximately equal to $r\log\frac{1}{r^2\lambda}$. We conclude that the minimum of (\ref{tominimize}) is asymptotically less then
$$\max\left(\sqrt{\log\frac{1}{\lambda r}}, r\log\frac{1}{\lambda r^2}\right).$$

%That means, that if $r(y)<R$ for all $y$ from the boundary of the polytope $P_1$ and for some fixed number $R$, then

We fix a positive number $R$ (which we will select later). Consider a convex polytope $P_1$ such that all its facets are close enough to the origin. In other words, assume that $r(y)<R$ for all $y\in \partial P_1$. Then

\begin{equation}\label{startfinal}
\gamma(\partial P_1)\lesssim \frac{\sqrt{n}}{t_0}\max_{r\in [0,R]}\left(\sqrt{\log\frac{1}{\lambda r}}, r\log\frac{1}{\lambda r^2}\right).
\end{equation}

We note that for $R\in (0,1)$, the right hand side of (\ref{startfinal}) is infinitely large. But as long as we assume that $R\in (1,\frac{1}{e\sqrt{\lambda}})$ the right hand side of (\ref{startfinal}) is asymptotically equal to

\begin{equation}\label{victory}
\frac{\sqrt{n}}{t_0}\max_{r\in [0,R]} \left(r\log\frac{1}{\lambda r^2}\right)\approx\frac{\sqrt{n}}{t_0}R\log\frac{1}{\lambda R^2},
\end{equation}

since $r\log\frac{1}{\lambda r^2}$ is increasing on $(1,\frac{1}{e\sqrt{\lambda}})$.

The estimate (\ref{victory}) is the first key ingredient for our proof. The other key ingredient is the following Lemma.

\begin{lemma}\label{hyper}
$$\gamma(H_{\rho})\lesssim\frac{\sqrt{n}}{t_0}\left(e^{-n}+ e^{-c\frac{n\rho^2}{t_0^2}}\right),$$
where $c$ is an absolute constant.
\end{lemma}
\textbf{Proof.} We write
$$\gamma(H_{\rho})=\frac{1}{n\nu_{n}J_{n-1}}\int_{\R^{n-1}} e^{-\varphi(\sqrt{|y|^2+\rho^2})}d\sigma(y).$$
Passing to the polar coordinates in $\R^{n-1}$, we get:
$$\gamma(H_{\rho})=\frac{(n-1)\nu_{n-1}}{n\nu_{n}J_{n-1}}\int_{0}^{\infty} s^{n-2}e^{-\varphi(\sqrt{s^2+\rho^2})}ds.$$
We make a change of variables $t=\sqrt{s^2+\rho^2}$ and use (\ref{balls_ratio}):
$$\gamma(H_{\rho})=\frac{(n-1)\nu_{n-1}}{n\nu_{n}J_{n-1}}\int_{\rho}^{\infty} t^{n-2}(1-\frac{\rho^2}{t^2})^{\frac{n-2}{2}}e^{-\varphi(t)}\frac{t}{\sqrt{t^2-\rho^2}}dt\approx$$
\begin{equation}\label{torefer}
\frac{\sqrt{n}}{J_{n-1}}\int_{\rho}^{\infty} t^{n-2}(1-\frac{\rho^2}{t^2})^{\frac{n-3}{2}}e^{-\varphi(t)}dt.
\end{equation}
It was shown in  \cite{KM} (Lemma 2.1) that
\begin{equation}\label{hvost}
\int_{5t_0}^{\infty} t^{n-2}e^{-\varphi(t)}dt\leq e^{-n} J_{n-2}.
\end{equation}
We note that $(1-\frac{\rho^2}{t^2})^{\frac{n-3}{2}}\leq 1$ for $n\geq 3$. Applying (\ref{hvost}) together with Lemma \ref{moments}, we conclude that for $n\geq 3$, (\ref{torefer}) is asymptotically smaller then
$$\frac{\sqrt{n}}{t_0}\left(e^{-n}+\max_{t\in [\rho,5t_0]} (1-\frac{\rho^2}{t^2})^{\frac{n-2}{2}}\right)\lesssim\frac{\sqrt{n}}{t_0}\left(e^{-n}+ e^{-c\frac{n\rho^2}{t_0^2}}\right).$$

For $n=2$ the surface area of any convex set is bounded by a constant. Thus for $n=2$ the result follows with the proper choice of $C$ in Theorem \ref{1}. This concludes the proof of the Lemma. $\square$
%since we assume that $\rho\leq \frac{t_0}{5\sqrt{\lambda n}}$.

Consider a polytope $P_2$ with $K$ facets such that all its facets are far enough from the origin. Namely, assume that $r(y)\geq R$ for all $y\in \partial P_2$. Then Lemma \ref{hyper} implies that

\begin{equation}\label{farfaces}
\gamma(\partial P_2)\lesssim \frac{\sqrt{n}}{t_0} K e^{-cR^2},
\end{equation}
as long as we chose $R\lesssim \sqrt{n}$.

Now we glue everything together. Let $R\in (1,\frac{1}{e\sqrt{\lambda}})$ (note that (\ref{farfaces}) is applicable for this range of $R$ since by Remark \ref{lambdas}$, \frac{1}{e\sqrt{\lambda}}\lesssim \sqrt{n}$). We split the surface of our polytope $P$ into two parts $P_1$ and $P_2$, where $P_1$ consists of the facets which are closer then $R$ to the origin and $P_2$ is the rest, i.e. the facets which are farther then $R$ from the origin. In other words,
$$P_1=\{y\in \partial P\,|\, r(y)\leq R\}$$
and
$$P_2=\{y\in \partial P\,|\, r(y)> R\}.$$
Applying (\ref{victory}) and (\ref{farfaces}) we observe, that
\begin{equation}\label{itall}
\gamma(\partial P)\lesssim \frac{\sqrt{n}}{t_0}\left(R\log\frac{1}{\lambda R^2}+ K e^{-R^2}\right).
\end{equation}
The estimate (\ref{itall}) holds for every $R\in (1,\frac{1}{e\sqrt{\lambda}})$. Minimizing (\ref{itall}) in $R$ we get that
\begin{equation}\label{itall1}
\gamma(\partial P)\lesssim \frac{\sqrt{n}}{t_0}\sqrt{\log K}\log\frac{1}{\lambda \log K}.
\end{equation}

Here we plugged $R\approx\sqrt{\log K}$, so the above estimate is valid for all $K\in [1, e^{\frac{c}{\lambda}}]$ for some absolute constant $c$. This finishes the proof of Theorem \ref{3}, and thus Theorem \ref{1}.$\square$

\section{Proof of the lower bound part}

Fix an integer $K\leq e^{\frac{c}{\lambda}}$. We consider K independent uniformly distributed random vectors $x_i\in \mathbb{S}^{n-1}$. Let $\rho\in (0, c\frac{t_0}{\sqrt{\lambda n}})$ (we will chose it later).
Consider a random polytope $P$ in $\R^{n}$, circumscribed around the ball of radius $\rho$:
$$P=\{x\in \R^{n} :\, \langle x,x_i\rangle\leq \rho, \,\,\,\forall i=1,...,K\}.$$
Passing to the polar coordinates as in Lemma \ref{hyper} and restricting the integration to $[t_0(1-\lambda), t_0(1+\lambda)]$ we estimate the expectation of $\gamma(\partial P)$ from below:
$$\mathbb{E}(\gamma(\partial P))\succsim$$
$$\frac{1}{n\nu_{n} J_{n-1}} K (n-1)\nu_{n-1} \int_{t_0(1-\lambda)}^{t_0(1+\lambda)} t^{n-2}e^{-\varphi(t)} (1-\frac{\rho^2}{t_0^2})^{\frac{n-3}{2}}(1-p(t))^{K-1}dt\succsim$$
\begin{equation}\label{expectation}
\frac{\sqrt{n}}{J_{n-1}} K \left(1-\frac{\rho^2}{t_0^2(1-\lambda)^2}\right)^{\frac{n-3}{2}}\int_{t_0(1-\lambda)}^{t_0(1+\lambda)} e^{-\varphi(t)} t^{n-2}(1-p(t))^{K-1}dt,
\end{equation}
where $p(t)$ is the probability that the fixed point on the sphere of radius $t$ is separated from the origin by the hyperplane $H_i=\{x:\,\langle x,x_i\rangle=\rho\}$. It was shown in \cite{GL}, Equation (70) (see also \cite{fedia}), that for $t\in [(1-\lambda)t_0,(1+\lambda)t_0]$
\begin{equation}\label{probability}
p(t)\precsim \frac{t_0}{\sqrt{n}\rho} \left(1-\frac{\rho^2}{t_0^2(1+\lambda)^2}\right)^{\frac{n-3}{2}}.
\end{equation}
We chose $\rho$ so that
$$K^{-1}=\frac{t_0}{\sqrt{n}\rho} \left(1-\frac{\rho^2}{t_0^2(1+\lambda)^2}\right)^{\frac{n-3}{2}},$$
which in the current range of $\rho$ means that $\rho=c\frac{t_0}{\sqrt{n}}\sqrt{\log K}$, and 
$$K\left(1-\frac{\rho^2}{t_0^2(1-\lambda)^2}\right)^{\frac{n-3}{2}}= \frac{\sqrt{n}\rho}{t_0}.$$
We use the above together with (\ref{expectation}) to conclude that the expectation $\mathbb{E}(\gamma(\partial P))$ is greater than
$$\frac{\sqrt{n}}{t_0}K(1-\frac{\rho^2}{t_0^2(1-\lambda)^2})^{\frac{n-3}{2}}\approx \frac{\sqrt{n}}{t_0}\frac{\sqrt{n}\rho}{t_0}=\frac{\sqrt{n}}{t_0}\sqrt{\log K},$$
which finishes the proof of Theorem \ref{4} and thus Theorem \ref{2}. $\square$

\section{Improvements in some partial cases}

In certain cases Theorem \ref{1} may be improved and made prescize. Namely, we fix $p>0$ and consider $\varphi(y)=\varphi_p(y)=\frac{|y|^p}{p}$ which corresponds to a measure $\gamma_p$ with density $e^{-\frac{|y|^p}{p}}$. Such measures are log-concave for $p\geq 1$. They were considered in \cite{GL1}. It was shown there, that for every convex body $Q$ in $\R^n$,
$$\gamma_p(\partial Q)\lesssim n^{\frac{3}{4}-\frac{1}{p}}.$$
The definition of $t_0$ implies that $t_0\approx n^{\frac{1}{p}}$ for the measures $\gamma_p$. Thus the above estimate can be rewritten:
$$\gamma_p(\partial Q)\lesssim \frac{\sqrt{n}}{t_0}n^{\frac{1}{4}}.$$
In particular, it was shown in \cite{GL1} that for every convex body $Q,$
\begin{equation}\label{sys2gl1}
\gamma_p(\partial Q)\lesssim \max_{y\in\partial Q}\frac{\alpha(y) |y|^{\frac{p}{2}}+1}{|y|^{\frac{p}{2}-1}}\approx \max_{y\in \partial Q}\frac{\sqrt{n}}{t_0}(r(y)+1),
\end{equation}

where, as before, $\alpha(y)=\cos(y,n_y)$ and $r(y)=\frac{\sqrt{n}}{t_0}\alpha(y)|y|.$ Using the scheme from the proof of Theorem \ref{3} we observe that for any polytope $P$ with $K$ facets, and for any $R>0,$
\begin{equation}\label{perfect}
\gamma_p(\partial Q)\lesssim \frac{\sqrt{n}}{t_0}(R+K e^{-cR^2}).
\end{equation}
Minimizing (\ref{perfect}) in $R,$ we get the following

\begin{theorem}\label{final}
$$\gamma_p(\partial P)\lesssim \frac{\sqrt{n}}{t_0}\sqrt{\log K}\approx n^{\frac{1}{2}-\frac{1}{p}}\sqrt{\log K}.$$
\end{theorem}

The above estimate is optimal since it coincides with the lower bound from Theorems \ref{2} and \ref{4}.

\end{document}